\documentclass[12pt]{article}

\usepackage{amsmath}
\usepackage{amsfonts}
\usepackage[T1]{fontenc} 
\newtheorem{theo}{Theorem}[section]
\newtheorem{prop}[theo]{Proposition}%[section]
\newtheorem{df}[theo]{Definition}%[section]
\newtheorem{lemme}[theo]{Lemma}%[section]
\newtheorem{Rem}[theo]{Remark}%[section]

\newcommand{\CQFD}{\hfill $\square$}

\newcommand{\ind}{\mathbf{1}}

\def\real{\Bbb{R}}

\def\ee{\Bbb{E}}

\def\E{\mathop{\hbox{\rm I\kern-0.20em E}}\nolimits}

\def\lto{-\!\!\!-\!\!\!\to}

\def\og{\leavevmode\raise.3ex
     \hbox{$\scriptscriptstyle\langle\!\langle$~}}
\def\fg{\leavevmode\raise.3ex
     \hbox{~$\!\scriptscriptstyle\,\rangle\!\rangle$}~}
     
\begin{document}

\title{Rescaled  weighted random balls models and stable self-similar random fields\footnote{This is a version of an original paper to be published in Stochastic Processes and Their Applications which only differs from the published paper by typographical changes (doi:10.1016/j.spa.2009.06.010).}}
\author{ Jean-Christophe Breton\footnote{Laboratoire MIA, Universit\'e de La Rochelle, 17042 La Rochelle Cedex, France. Email: jcbreton@univ-lr.fr} 
\ and Clément Dombry\footnote{Laboratoire LMA, Universit\'e de Poitiers, 
 Téléport 2, BP 30179, F-86962 Futuroscope-Chasseneuil cedex, France. 
 Email: clement.dombry@math.univ-poitiers.fr}
 \protect\hspace{1cm}}
\date{}
\maketitle

\abstract{ 
We consider weighted random balls in $\real^d$ distributed according to a random Poisson measure with heavy tailed intensity and study the asymptotic behavior of the total weight of some configurations in $\real^d$ while we perform a zooming operation. 
The resulting procedure is very rich and several regimes appear in the limit, depending on the intensity of the balls, the zooming factor, the tail parameters of the radii and the weights. 
Statistical properties of the limit fields are also evidenced, such as isotropy, self-similarity or dependence. 
One regime is of particular interest and yields  $\alpha$-stable stationary isotropic self-similar generalized random fields which recovers Takenaka fields, Telecom process or fractional Brownian motion.
%The asymptotic behaviors of the total weight of some configurations in $\real^d$ are studied. 
%Different scaling, yielding different limit fields, are investigated.
} 
\\

\section*{Introduction}
%\subsection{Motivations}

In this work, we consider the so-called weighted random balls model and investigate its convergence when suitably rescaled and normalized. 
We exhibit  three different asymptotic regimes driving the macroscopic and microscopic variations of this model, 
namely  
$(i)$ a stable, translation and rotation invariant, self-similar random field on $\real^d$, 
$(ii)$ a Poissonian field 
and $(iii)$ a stable field with independence. 
The weighted random balls model is constructed in the following way: the centers of the balls are distributed according to a Poisson point process, with each center $x$ labelled with a random radius $r$ and a random weight $m$. 
The field under study is, roughly speaking,
at each point, the weight density defined as the sum of the weights of the balls containing this point. 
The overlap of the balls yields non-trivial spatial correlations when the random radii of the balls are heavy tailed.

This fairly simple geometric construction has found numerous applications and is pertinent in various modeling situations. 
Similar stochastic models were considered by Kaj in \cite{Kaj-net} when modeling a simplified wireless network that consists of a collection of spatially distributed stations equipped with emitters  for transmission over a common communication channel. 
Here, the location of a station or of a network node is represented by the point $x$, its range by the radius $r$ and its power by the weight $m$. 
The weight density measures the total power of emission at a given point
and in this case, $m$ is supposed to be non-negative. But our model supports more generally real-valued weights. 

In \cite{BE}, Biermé and Estrade consider similar models in dimension $d=2$ as models in imagery (in this case, the weight intensity stands for the gray level of a pixel in a black and white picture)
and in dimension $d=3$ for modeling three dimensional porous or  heterogeneous media (here, the weight density is seen as a mass density). 
They investigate the microscopic properties of the random balls configurations 
by performing a scaling operation which amounts to zoom in smaller regions of space. 
In \cite{KLNS}, Kaj {\it et al.} study similar random grain model by shrinking to zero the volume of the grains. 
This amounts to analyse the macroscopic properties of the random balls configurations by 
performing a scaling operation which amounts here to zoom out over larger areas. 

Recently, Biermé, Estrade and Kaj introduce in \cite{BEK} a general framework for rescaled random balls model
allowing both zoom-in (as in \cite{BE}) and zoom-out (as in \cite{KLNS}). 
In this zooming procedure, several limit fields arise, which are either of Gaussian or of Poisson type according to the  respective asymptotic of the zooming rate and of the Poisson intensity of the balls. 
Furthermore, they show that essentially all Gaussian, translation and rotation invariant self-similar generalized random fields can be obtained as such a limit.

Note that in the rescaled random balls model of \cite{BE}, \cite{KLNS} and \cite{BEK}, 
the weights in the field under study are fixed equal to $m\equiv 1$. 
% exhibiting several scaling limits giving non-trivial limit fields.
Models with randomized weights have been less intensively studied. 
In dimension $d=1$, Kaj and Taqqu  study in \cite{KT} limiting schemes for weighted random balls model, deriving Gaussian, Poisson and stable regimes. 
This model applies in particular to study the random variation in packet networks computer traffic.

Our main contribution in this paper is to introduce a general study of macroscopic and microscopic variations in weighted models in $\real^d$. 
This generalizes both \cite{BEK} since the balls are randomly weighted 
and \cite{KT} since we consider an arbitrary dimension $d$ and more general configurations on the balls. 
As in \cite{KLNS} and \cite{KT}, three different regimes appear according to 
the relative behavior of the scaling rate and of the Poisson intensity. 
In particular, when the random weights are heavy tailed, the limit generalized random fields are stable, translation and rotation invariant, and also self-similar.
The paper is organized as follows. 
The model under study is described in Section~\ref{sec:model}. 
Our main results under different scaling regimes are stated and discussed in Section~\ref{sec:result}. 
Finally, Section~\ref{sec:proof} is devoted to the proof of technical lemmas and of the main results. 

%%%%%%%%%%%%%%%%%%%%%%%%%%%%%%%%%%%%%%%%%%%%%%%%%%%%%%%%%%%%%%%%%%%%%%%%%%%%%%%%%%%%%%%%%%%%%%%%%%%%%%%%%%%%%%%%%%

\section{Model of weighted random balls}
\label{sec:model}
We consider random balls $B(x,r)=\{y\in \real^d\::\: \|y-x\|<r\}$ with weight $m$, 
the triplet $(x,r,m)$ being distributed according to a Poisson random measure $N_\lambda(dx,dr,dm)$ on $\real^d\times\real^+\times\real$ with intensity 
$$
n(dx,dr,dm)=\lambda dx F(dr) G(dm)
$$
where $\lambda$ is positive, 
$F$ is a positive measure on  $\real^+$ and $G$ a probability measure on  $\real$. 
Here, and in what follows, $\|\cdot\|$ stands for the usual Euclidean norm on $\real^d$.

%M redaction pour \lambda
The point process of the centers of the balls in $\real^d$ is the projection of the point process in $\real^d\times\real^+\times\real$ corresponding to the Poisson random measure $N_\lambda(dx,dr,dm)$. 
It is easily seen that it is a Poisson point process with intensity 
$\lambda dx$, and hence the parameter $\lambda$ is interpreted as the intensity of the balls in $\real^d$.

%M redaction pour F
 We suppose that the  measure $F$ driving the distribution of the radius $r$ is absolutely continuous  $F(dr)=f(r)dr$ with 
\begin{equation}
\label{eq:F}
 \int_{\real^+}r^dF(dr)<+\infty
\end{equation}
and such that for either $\epsilon=+ 1$ or $\epsilon=-1$, 
\begin{equation}
\label{eq:F1}
f(r)\sim_{r\to 0^\epsilon} C_\beta r^{-1-\beta} 
\end{equation}
where by convention $0^{+1}=0$ and $0^{-1}=+\infty$. As will be explained later, the case $\epsilon=+1$ will be referred as the zoom-in case, whereas the case $\epsilon=-1$ will be referred as the zoom-out case. 
Condition \eqref{eq:F1} assumes a power behavior of the radius density at the origin (zoom-in case $\epsilon=+1$) or at infinity (zoom-out case $\epsilon=-1$). 
Condition \eqref{eq:F} is equivalent to the finiteness of the volume of the random balls.
Note that  assumptions \eqref{eq:F} and \eqref{eq:F1} together imply that for $\epsilon=+1$, we must have $\beta<d$, while for $\epsilon=-1$, we must have $\beta>d$.

%M redaction pour G
We suppose that the probability measure $G$ belongs to the normal domain of attraction of the $\alpha$-stable distribution $S_\alpha(\sigma,b,\tau)$ with $\alpha\in (1,2]$, 
i.e. if $X_1, \dots, X_n$ are independent and identically distributed (i.i.d.) according to $G$, $n^{-1/\alpha}(X_1+\dots+X_n)\Rightarrow  S_\alpha(\sigma,b,\tau)$. 
We recall the following estimate (see \cite{Feller}) of the characteristic function $\varphi_G$ of $G$ as $\theta\to 0$ 
\begin{equation}
 \varphi_G(\theta)=1+i\theta\tau-\sigma^\alpha|\theta|^\alpha(1+i b \varepsilon(\theta)\tan(\pi\alpha/2)+o(|\theta|^\alpha),
\end{equation}
where here, and in what follows, $\varepsilon(a)=+1$ if $a>0$, $\varepsilon(a)=-1$ if $a< 0$ and $\varepsilon(0)=0$. 
In case $\alpha\in (1,2)$, typical choices for $G$ are heavy tailed distributions 
while for $\alpha=2$, $G$ is any distribution with finite variance.
In this latter case, we recover a weighted version of the main results in \cite{BEK} 
(set $G=\delta_1$ to recover exactely the setting described in \cite{BEK}).

\medskip
\noindent
Let $\cal{M}$ denote the set of signed measures on $\mathbb{R}^d$ with finite total variation $|\mu|(\mathbb{R}^d)$, where $|\mu|$ is the total variation of a measure $\mu$. We recall that equipped with the norm of total variation $\|\mu\|_{\cal M}=|\mu|(\mathbb R^d)$,  $\cal{M}$ is a Banach space. %({\it cf.} \cite{Rudin} p. 161). 
We consider the random field
\begin{equation}
\label{eq:M}
M(\mu)=\int_{\real^d\times\real^+\times\real} m \mu(B(x,r)) N_\lambda(dx,dr,dm)
\end{equation}
indexed by signed measures $\mu\in\cal{M}$. 
When $\mu=\delta_y$, $y\in\real^d$, $M(\delta_y)$ is the weight density at point $y$ as described in the introduction: it is the sum of the algebraic weights of the balls containing the point $y$. 

Note that the stochastic integral in \eqref{eq:M} is well defined and has finite expected value since
\begin{eqnarray*}
& &\int_{\real^d\times\real^+\times\real} \hskip-25pt|m \mu(B(x,r))| n(dx,dr,dm) \\
&\leq& \int_{\real} |m| G(dm) \times \lambda |B(0,1)| |\mu|(\real^d) \int_{\real^+} r^d F(dr) <+\infty.
\end{eqnarray*}
where $|A|$ stands for the Lebesgue measure of a Borel set $A$. 
Furthermore, the expected value is given by
$$
\ee[M(\mu)]= \lambda |B(0,1)| \int_{\real} m G(dm)  \int_{\real^+} r^d F(dr)\ \mu(\real^d).
$$
%The random field $M$ induces a linear functional ${\cal M}\rightarrow L^1( \Omega,{\cal F},\mathbb{P})$, $ \mu\mapsto M(\mu)$. 
%Moreover
%$$
%\|M(\mu)\|_{L^1}=\ee[|M(\mu)|]\leq \left(\lambda |B(0,1)| \int_{\real} |m| G(dm)  \int_{\real^+} r^d F(dr)\right) \|\mu\|_{\cal M}.
%$$
%This shows that the linear functional is bounded and gives an upper bound of its operator norm.

\medskip
We are interested in the variations of $M(\mu)$ at a microscopic or macroscopic level. 
To do so, we swell, resp. shrink, the volume of the balls replacing the radius $r$ of a ball by $\rho r$ and taking the limit $\rho\to +\infty$, resp. $\rho\to 0$. 
In this procedure, the law of the radius is replaced by $F_\rho(dr)=f(r/\rho) dr/\rho$, the image measure of $F(dr)$ by the change of scale $r\mapsto \rho r$.  
In order to derive non-trivial asymptotics, the intensity $\lambda$ of the balls is changed accordingly and we shall write $\lambda(\rho)$ 
to underline that from now on the intensity depends on the scaling parameter $\rho$. 
In what follows, we are thus interested in the following random field:
$$
M_{\rho}(\mu)=\int_{\real^d\times \real^+\times\real} m \mu(B(x,r)) N_{\lambda(\rho),\rho}(dx,dr,dm)
$$ 
where $N_{\lambda(\rho),\rho}(dx,dr,dm)$ is the Poisson measure with intensity $\lambda(\rho)dx F_\rho(dr)G(dm)$. 
The limit $\rho\rightarrow 0$ is interpreted as zoom-out in the random configurations of balls 
and this is relevant when the behavior of $f$ is known at $+\infty$, i.e. $\epsilon=-1$ in \eqref{eq:F1}.
In this case, we investigate the macroscopic variations of $M$.
On the contrary, $\rho\rightarrow +\infty$ is interpreted as zoom-in in space and this is 
relevant when the behavior of $f$ is known at $0$, i.e. $\epsilon=+1$ in \eqref{eq:F1} and this is the microscopic variations that are investigated.

\begin{Rem}
\label{rem:KT1}
{\rm As observed before, the choice $G=\delta_1$ recovers the setting of \cite{BEK} for non-weighted random balls, see (4) therein. 
If $d=1$, a {\it verbatim} replacement of $B(x,r)=(x-r,x+r)$ by $(x, x+r)$ and the choice $\mu= |\cdot \cap (0,t)|$ recover the field studied in \cite{KT} in the "continuous flow reward model", see (18) therein. 
}
\end{Rem}

%%%%%%%%%%%%%%%%%%%%%%%%%%%%%%%%%%%%%%%%%%%%%%%%%%%%%%%%%%%%%%%%%%%%%%%%%%%%%%%%%%%%%%%%%%%%%%%%%%%%%%%%%%%%%%%%%%

\section{Results}
\label{sec:result}

We exhibit normalization terms $n(\rho)$ such that the normalized centered random field $n(\rho)^{-1}\big(M_\rho(\cdot)-\ee[M_\rho(\cdot)]\big)$ converges in finite-dimensional distribution ({\it f.d.d.}) to a limit random field. 
In what follows, we are interested in {\it f.d.d.} convergence on subspaces $\widetilde {\cal M}$ of ${\cal M}$ and we will denote it by $\stackrel{\widetilde {\cal M}}{\longrightarrow}$. 

It is natural to investigate first the behavior of the random field giving the density of the weights at each point which in our notations rewrites $(M_\rho(\delta_y))_{y\in\real^d}$. The heuristic is the following. The average numbers of balls containing the point $y$ is given by 
$$
\ee\left[\int_{\real^d\times\real^+\times\real} \ind_{\{y\in B(x,r)\} }N_{\lambda(\rho),\rho}(dx,dr,dm)\right]=V \lambda(\rho)\rho^d,
$$
where $V=c_d \int r^d F(dr)$ is the expected volume of a random ball
and $c_d$ stands for the volume of the Euclidean unit ball in $\real^d$. 
Since the weights belong to the domain of attraction of an $\alpha$-stable distribution, it is natural to introduce the scaling $n_0(\rho)= \lambda(\rho)^{1/\alpha}\rho^{d/\alpha}$. Convergence of the normalized and centered random variable $M_\rho(\delta_y)$ to an $\alpha$-stable distribution is obtained if we suppose  that  $\lambda(\rho)\rho^d\to +\infty$ when $\rho\to 0 ^{-\varepsilon}$. Heuristically, the dependence between $M_\rho(\delta_{y_1})$ and $M_\rho(\delta_{y_2})$ is given by the weights of the balls containing both points $y_1$ and $y_2$.  In the zoom-in case ($\varepsilon=-1, \rho\to+\infty$), the balls are very large yielding total dependence at the limit and we have:
\begin{equation}\label{eq:triv1}
n_0(\rho)^{-1}(X_\rho(\delta_y)-\ee[X_\rho(\delta_y)]) \stackrel{f.d.d.}{\lto}  W_\alpha, \quad y\in\real^d
\end{equation}
where $W_\alpha(y)\equiv W_\alpha$ is a constant random field distributed according to $S_\alpha(\sigma V^{1/\alpha},b,0)$.
In the zoom-out case ($\varepsilon=-1, \rho\to 0$), the balls are very small yielding independence at the limit and we have:
\begin{equation}\label{eq:triv2}
n_0(\rho)^{-1}(M_\rho(\delta_y)-\ee[M_\rho(\delta_y)]) \stackrel{ f.d.d.}{\lto} W_\alpha(\delta_y), \quad y\in\real^d,
\end{equation}
where $W_\alpha(\delta_y), y\in\real^d$, are i.i.d.  $S_\alpha(\sigma V^{1/\alpha},b,0)$ distributed. 
Similar results as in \eqref{eq:triv1} and in \eqref{eq:triv2} hold true for f.d.d. convergence on the space of measures with finite support. 
Since these results are not surprising, their proofs are omitted and in what follows we investigate convergence for more general measures.

%M Supprimé une grande partie
%We can show that when $n_0(\rho):=\lambda(\rho)\rho^d\to+\infty$, for all $y\in \real^d$:
%\begin{equation}
%\label{eq:limit_dirac}
%n_0(\rho)^{-1}\big(M_\rho(\delta_y)-\ee[M_\rho(\delta_y)]\big)\Rightarrow 
%S_\alpha\left(\sigma c_d^{1/\alpha}\left(\int_{\real^d} r^d F(dr)\right), b,0\right)
%\end{equation} 
%where $\sigma$ is given in \eqref{eq:sigma}, $b$ in \eqref{eq:biais} ????
%and $c_d$ stands for the volume of the unit ball in $\real^d$ and is explicitely given by 
%$$
%c_d=\left\{
%\begin{array}{ll}
%\frac{2^{2p+1}\pi^p p!}{(2p+1)!} & \mbox{when } d=2p+1,\\
%\frac{\pi^p}{p!} & \mbox{when } d=2p.
%\end{array}
%\right .
%$$
%Since the limit does not depend on $y\in\real^d$, we do not have {\it f.d.d.} convergence for 
%$n(\rho)^{-1}\big(M_\rho(\delta_y)-\ee[M_\rho(\delta_y)]\big)$. 
%Furthermore as a consequence of \eqref{eq:limit_dirac}, the following proposition gives an 
%indication that no limit theorem can be given on the whole space ${\cal M}$ at least for smaller normalization than $n_0(\rho)$. 

%M je supprime Banach-Steinhaus : ca souleve plus de questions que cela n'en resoud.

%\begin{prop}
%\label{Prop:Banach-Steinhaus}
%Let $n(\rho)=o(n_0(\rho))$. 
%There is a dense subset of $\cal M$ on which the rescaled process 
%$$
%n(\rho)^{-1}\big(M_\rho(\mu)-\ee[M_\rho(\mu)]\big) 
%$$ 
%diverges in $L^1(\Omega,\cal F,\mathbb P)$ as $\rho\rightarrow 0^{-\epsilon}$.
%\end{prop}

%M Nouvelle transition

%%%%%%%%%%%%%%%%%%%%%%%%%%%%%%%%%%%%%%%%%%%%%%%%%%%%%%%%%%%%%%%%%%%%%%%%%%%%%%%%%%%%%%%%%%%%%%%%%%%%%%%%%%%%%%%%%%

\subsection{Preliminaries on measured spaces}
We introduce a subspace  ${\cal M}_{\alpha,\beta}\subset {\cal M}$ on which we will show the convergence of the rescaled generalized random field $M_\rho(\mu)$.
\begin{df} For $1<\alpha\leq 2$ and $\beta>0$, let ${\cal M}_{\alpha,\beta}$ be the subset of measures $\mu\in {\cal M}$ satisfying for some finite constant $C$ and some $0<p<\beta<q$:
\begin{equation}
\label{eq:Cmu}
\gamma(r):=\int_{\real^d}|\mu(B(x,r))|^\alpha dx \leq C(r^p\wedge r^q) %\quad \mbox{ for some } p<\beta<q. 
\end{equation}
where for reals $a, b$: $a\wedge b=\min(a,b)$. 
\end{df}

Here, and in what follows, $C$ is a finite constant that may change at each occurrence.
Some elementary properties of the spaces ${\cal M}_{\alpha,\beta}$ are given in the following proposition.
\begin{prop}
\label{prop:subspace1} \  
\begin{description}
\item{i)} ${\cal M}_{\alpha,\beta}$ is a linear subspace of ${\cal M}$ on which 
$$\forall \mu\in {\cal M}_{\alpha,\beta}\ ,\ \int_{\real^d\times\real^+}|\mu(B(x,r))|^\alpha r^{-\beta-1}dxdr <+\infty.$$ 
\item{ii)} ${\cal M}_{\alpha,\beta}$ is closed under translations, rotations and dilatations, {\it i.e.}
when $\mu\in{\cal M}_{\alpha,\beta}$, $\tau_s\mu$, $\Theta\mu$ and $\mu_a$ are also in ${\cal M}_{\alpha,\beta}$
where for any Borelian set $A$ and for $s\in\real^d$, $\Theta\in {\cal O}(\real^d)$, $a\in\real_+$ 
$$
\tau_s\mu(A)=\mu(A-s), 
\Theta\mu(A)=\mu(\Theta^{-1}A), 
\mu_a(A)=\mu(a^{-1}A).
$$
%\item{iii)} If $d-1<\beta<\alpha d$, then $\delta_y-\delta_x\in {\cal M}_{\alpha,\beta}$ for any $(x,y)\in\real^d$.
\item{iii)} When $\alpha\leq \alpha'$, we have ${\cal M}_{\alpha,\beta}\subset{\cal M}_{\alpha',\beta}$.
\item{iv)} When $\beta\geq d$, the space ${\cal M}_{\alpha,\beta}$ is included in the subspace of diffuse measures (i.e. such that $\mu(\{x\})=0$ for any $x\in\real^d$).
\item{v)} When $\beta\leq d$, the space ${\cal M}_{\alpha,\beta}$ is included in the subspace of centered measures (i.e. such that $\mu(\real^d)=0$).
\end{description}
\end{prop}
Observe that Dirac measures $\delta_y$, $y\in \real^d$,  are not in ${\cal M}_{\alpha,\beta}$. 
However, explicit examples of measure in ${\cal M}_{\alpha,\beta}$ are given in the following proposition. 
Absolutely continuous measures (with respect to the Lebesgue measure) $\mu(dx)=\phi(x)dx$ with integrable density $\phi\in L^1(\real^d)\cap L^\alpha(\real^d)$ will play an important role. In this case, we shall (abusively) note $\mu\in L^1(\real^d)\cap L^\alpha(\real^d)$.
\begin{prop}
\label{prop:subspace2} \  
\begin{description}
\item{i)} If $d<\beta<\alpha d$,  any measure $\mu\in L^1(\real^d)\cap L^\alpha(\real^d)$ belongs to ${\cal M}_{\alpha,\beta}$.
%for $\eta>d^2/\beta$  %sinon la densite n'est pas integrable.
\item{ii)} If $d-1<\beta < d$, any centered measure $\mu(dx)=\phi(x)dx\in L^1(\real^d)$ such that $\int_{\real^d}|\!|y|\!||\phi(y)|dy <+\infty$ belongs to ${\cal M}_{\alpha,\beta}$, as well as any centered measure with finite support.
%for $\eta>d^2/\beta$  %sinon la densite n'est pas integrable.
\end{description}
\end{prop}
Note that in particular, when $d<\beta<\alpha d$ (resp. $d-1<\beta<d$), ${\cal M}_{\alpha,\beta}$ contains the space $\cal S$ of measures with density in the Schwartz class  (resp. ${\cal S}_0$ the space of centered measures with density in the Schwartz class). 
Note also that when $\alpha=2$, the conditions supposed in \cite{BEK} on the measure $\mu$ (expressed in terms of Riesz energy) imply that $\mu\in{\cal M}_{2,\beta}$. 
By analogy with the case $\alpha=2$, we suspect the space ${\cal M}_{\alpha,\beta}$ to be reduced to $\{0\}$ whenever $\beta\leq d-1$ or $\beta\geq \alpha d$, but we have no formal proof of these facts. However, we refer to Theorem \ref{theo:indep-stable2} for a positive result when $\beta>\alpha d$.

%%%%%%%%%%%%%%%%%%%%%%%%%%%%%%%%%%%%%%%%%%%%%%%%%%%%%%%%%%%%%%%%%%%%%%%%%%%%%%%%%%%%%%%%%%%%%%%%%%%%%%%%%%%%%%%%%%

\subsection{Limit theorems for the rescaled weighted random balls model}

We now come to the main results of this paper, {\it viz.} limit theorems for the rescaled generalized random fields $M_\rho$ and for configurations $\mu\in{\cal M}_{\alpha,\beta}$ on the balls.
As in \cite{KLNS} and \cite{KT} (for $\epsilon=-1$), 
several regimes appear according to the density of large/small balls in the limit. 
More precisely, using \eqref{eq:F1}: 

{\bf Zoom-out case} ($\epsilon=-1$, i.e. $\beta>d$ and $\rho\to 0$). The mean number of balls with radius larger than one that cover the origin is given by
$$
\int_{\real^d\times\real+}\ind_{\|x\|<r}\ind_{r>1}\lambda(\rho)dxF_\rho(dr)
=c_d\lambda(\rho)\int_1^{+\infty} r^dF_\rho(dr) \sim_{\rho\to 0} \frac{c_dC_\beta}{\beta-d} \lambda(\rho)\rho^\beta.
$$
Consequently, we distinguish the following three scaling regimes:
\begin {itemize}
\item large-balls scaling: $\lambda(\rho)\rho^\beta\to+\infty$, 
\item intermediate scaling: $\lambda(\rho)\rho^\beta\to a\in(0,+\infty)$,
\item small-balls scaling: $\lambda(\rho)\rho^\beta\to0$.
\end{itemize}

{\bf Zoom-in case} ($\epsilon=+1$, i.e. $\beta<d$ and $\rho\to +\infty$). The mean number of balls with radius less than one that cover the origin is given by
$$
\int_{\real^d\times\real+}\ind_{\|x\|<r}\ind_{r<1}\lambda(\rho)dxF_\rho(dr)=c_d\lambda(\rho)\int_0^1 r^dF_\rho(dr) \sim_{\rho\to +\infty} \frac{c_dC_\beta}{d-\beta} \lambda(\rho)\rho^\beta.
$$
In this case, the three scaling regimes are:
\begin{itemize}
\item small-balls scaling: $\lambda(\rho)\rho^\beta\to+\infty$, 
\item intermediate scaling: $\lambda(\rho)\rho^\beta\to a\in(0,+\infty)$,
\item large-balls scaling: $\lambda(\rho)\rho^\beta\to 0$.
\end{itemize}

\bigskip
In what follows, we study precisely the limiting shape of the random balls by investigating the fluctuations of $M(\mu)$ around its mean. 
Three different limit fields are exhibited according to the scaling performed. 

%%%%%%%%%%%%%%%%%%%%%%%%%%%%%%%%%%%%%%%%%%%%%%%%%%%%%%%%%%%%%%%%%%%%%%%%%%%%%%%%%%%%%%%%%%%%%%%%%%%%%%%%%%%%%%%%%%

\subsubsection{Stable regime with dependence}
\label{sec:stable1}
In this section, we investigate the behavior of $M$ under the scaling $\rho^\beta\lambda(\rho)\to +\infty$. 
In this case the limiting field is given by an $\alpha$-stable integral. 
We recall that the stable stochastic integral of $f$ with respect to an $\alpha$-stable random measure with control measure $m$ is well defined whenever $f\in L^\alpha(dm)$
and in this case, this stochastic integral follows an $\alpha$-stable distribution. 
We refer to \cite{ST} for a complete account on stable measures and integrals. 
The asymptotic of the rescaled generalized fields $M_\rho$ is given by the following result:
\begin{theo}
\label{theo:stable}
Suppose $\rho^\beta\lambda(\rho)\to +\infty$ when $\rho\to 0^{-\epsilon}$. 
Let $n_1(\rho)=\lambda(\rho)^{1/\alpha} \rho^{\beta/\alpha}$. 
We have 
\begin{equation}
\label{eq:stable}
\frac{M_\rho(\cdot) -\ee[M_\rho(\cdot)]}{n_1(\rho)} \stackrel{{\cal M}_{\alpha,\beta}}{\lto} Z_\alpha(\cdot)\quad \rho\to0^{-\epsilon}
\end{equation}
where $Z_\alpha(\mu)=\int_{\real^d\times\real^+} \mu(B(x,r))M_\alpha(dr,dx)$ is a stable integral with respect to the $\alpha$-stable measure $M_\alpha$ with control measure
$\sigma^\alpha C_\beta r^{-1-\beta} drdx$ and constant skewness function $b$ given in the domain of attraction of $G$. 
\end{theo}
Note that $Z_\alpha(\mu)$ makes sense as soon as 
$\int_{\real\times\real^+ }|\mu(B(x,r))|^\alpha r^{-1-\beta} drdx <+\infty$ (see Prop. \ref{prop:subspace1}-$i)$).
However, we need the stronger assumption $\mu\in{\cal M}_{\alpha,\beta}$ in order to derive \eqref{eq:stable}.
Roughly speaking, the control \eqref{eq:Cmu} of $\mu\in{\cal M}_{\alpha, \beta}$ allows to replace $F$ by its tails behavior given in \eqref{eq:F1} in asymptotic estimate. 

Due to the invariance by translation and rotation of the Lebesgue measure, the self-similarity of stable integral and the (global) invariance by rotation of the balls 
and because of Proposition \ref{prop:subspace1}-ii), we derive the following properties for the limit field $Z_\alpha$ of Theorem \ref{theo:stable}:
\begin{prop}\ 
\label{prop:properties_stable}
\begin{description}
\item{i)} The field $Z_\alpha$ is stationary on ${\cal M}_{\alpha,\beta}$, that is: 
$$
\forall \mu\in {\cal M}_{\alpha,\beta}, \forall s\in \real^d, Z_\alpha(\tau_s\mu)\stackrel{fdd}{=} Z_\alpha(\mu).
$$
\item{ii)} The field $Z_\alpha$ is isotropic on ${\cal M}_{\alpha,\beta}$, that is:
$$
\forall \mu\in {\cal M}_{\alpha,\beta}, \forall \Theta\in {\cal O}(\real^d), Z_\alpha(\Theta\mu)\stackrel{fdd}{=} Z_\alpha(\mu).
$$
\item{iii)} The field $Z_\alpha$ is self-similar on ${\cal M}_{\alpha,\beta}$ with index $(d-\beta)/\alpha$, 
that is: 
$$
\forall \mu\in {\cal M}_{\alpha,\beta}, \forall a>0, Z_\alpha(\mu_a)\stackrel{fdd}{=} a^{(d-\beta)/\alpha}Z_\alpha(\mu).
$$
\end{description}
\end{prop}

\begin{Rem}
{\rm
The covariation gives an insight into the structure of the spatial dependence of the stable generalized field. 
It is a generalization of the usual notion of covariance to the stable framework.  
Here, for $\mu_1,\mu_2\in {\cal M}_{\alpha,\beta}$, the covariation of $Z_\alpha(\mu_1)$ on $Z_\alpha(\mu_2)$ is given by
\begin{eqnarray*}
& &[Z_\alpha(\mu_1),Z_\alpha(\mu_2)]_\alpha\\
&=&\sigma^\alpha C_\beta\int_{\real^d\times\real^+} \mu_1(B(x,r))\epsilon(\mu_2(B(x,r)))|\mu_2(B(x,r))|^{\alpha-1}r^{-\beta-1}dr dx.
\end{eqnarray*}
Note that the integral above is well defined by Hölder's inequality since $\mu_1$ and $\mu_2$ belong to 
${\cal M}_{\alpha,\beta}$. 
We refer to \cite{ST} for a definition and properties of the covariation. 
Note that unlike the Gaussian case, the covariation structure is not sufficient to characterize the distribution of the generalized random field. 
However, since even if $\mu_1$ and $\mu_2$  have disjoint supports,  $[Z_\alpha(\mu_1),Z_\alpha(\mu_2)]_\alpha\neq 0$, $Z_\alpha(\mu_1)$ and $Z_\alpha(\mu_2)$ are not independent and the random field $Z_\alpha$ is stable with dependence. 
}

\end{Rem}

%\begin{Rem}
%{\rm
%We recall the following notion of long-range dependence (also called long-memory) 
%for one-dimensional heavy tailed processes from Yang-Petropulu \cite{YP}:
%an $\alpha$-stable stationnary process $X(t)$ is said to be long-range dependent
%if for some $c>0$ and $0<\gamma<1$.  
%$$\lim_{\tau\to\infty}|\tau|^{1-\gamma}[X(t),X(t+\tau)]_\alpha = c.$$
%In the present framework, for any $\mu\in{\cal M}_{\alpha,\beta}$, the process $Z_\alpha(\tau_s\mu)$, $s\in\real^d$, is %$\alpha$-stable and stationnary. 
%We suspect that for $\mu\neq 0$, this process exhibits long-range dependence. 
%Unfortunately, we only have a partial result in this direction: 
%when $\beta>d$, we can prove that if $\mu$ is non-negative and compactly supported, the covariation decreases slowly. More %precisely, let $\beta>d$ and let $\mu\in {\cal M}_{\alpha,\beta}$ be non-negative and compactly supported. Then,  
%$$
%\liminf_{s\to\infty}|s|^{\beta-d}[Z_\alpha(\mu),Z_\alpha(\tau_s\mu)]_\alpha >0.
%$$
%To see this, let $M$ be such that $B(0,M)$ contains the support of $\mu$. Note that if $r>M+|s|$ and $|x|\leq r-M-|s|$, 
%then $\mu(B(x,r))=\tau_s\mu(B(x,r))=\mu(\real^d)>0$, and hence
%\begin{eqnarray*}
%[Z_\alpha(\mu),Z_\alpha(\tau_s\mu)]_\alpha &\geq& \sigma^\alpha C_\beta \int_{r>M+|s|}\int_{|x|\leq r-M-|s|} %\mu(\real^d)^\alpha
%r^{-\beta-1}dxdr \\
%& =& \sigma^\alpha C_\beta \mu(\real^d)^\alpha c_d \sum_{k=0}^d (-1)^k {d\choose k}(k+\beta-d)^{-1} (M+|s|)^{d-\beta} \\
%&\geq & c|s|^{d-\beta}
%\end{eqnarray*}
%for large $|s|$ and some $c>0$.
%}
%\end{Rem}
 
\begin{Rem}
{\rm
Note that when $d-1<\beta<d$, $\mu_z=\delta_z-\delta_0$ for $z\in\real^d$ belongs to ${\cal M}_{\alpha,\beta}$. 
For such a measure, when moreover $b=0$ (i.e. when $G$ in our model is symmetric), our limiting field rewrites 
$$
Z_\alpha(\mu_z)=\int_{\real^d\times \real^+} \ind_{B(z,r)\Delta B(0,r)}M_{\alpha}(dx,dr)
$$
where $A\Delta B=(A\setminus B)\cup (B\setminus A)$.
In this case, we recover the so-called $(\alpha, H)$-{\it Takenaka field} with $H=(d-\beta)/\alpha$.   
It is self-similar with index $H$,  with stationary increments and almost surely with continuous sample paths, 
see \cite[p. 25]{BEK} or \cite[Sect. 8.4]{ST}. 
}
\end{Rem}

\begin{Rem}
\label{rem:KT_stable_dep}
{\rm
When $d=1$, $\beta\in(1,\alpha)$ and $\mu_t=|\cdot\cap(0,t)|$, 
the field $Z_\alpha(\mu_t)$ coincides with the {\it Telecom process} obtained in the fast connection rate for the "continuous flow reward model"
in \cite[Th. 2]{KT}, see also Remark \ref{rem:KT1} above. 
% Petite verification \mu_t est bien dans {\cal M}_{\alpha,\beta}: 
%Observe that $\mu_t \in {\cal M}_{\alpha,\beta}$. Indeed, for $r$ small, 
%\begin{eqnarray*}
%\gamma(r)&=&\int_{-r}^r|(0, x+r)|^\alpha dx+\int_r^{t-r}|(x-r,x+r)|^\alpha dx+\int_{t-r}^{t+r}|(x-r,t)|^\alpha dx\\
%&=&\frac{(2r)^{\alpha+1}}{\alpha+1}+(t-2r)(2r)^{\alpha}+\frac{(2r)^{\alpha+1}}{\alpha+1}
%=O(r^\alpha)
%\end{eqnarray*}
%and for $r$ large:
%\begin{eqnarray*}
%\gamma(r)&=&\int_{-r}^{t-r}|(0, x+r)|^\alpha dx+\int_{t-r}^rt^\alpha dx+\int_r^{t+r}|(x-r,t)|^\alpha dx\\
%&=&\frac{t^{\alpha+1}}{\alpha+1}+t^\alpha (2r-t)+\frac{t^{\alpha+1}}{\alpha+1}
%=O(r).
%\end{eqnarray*}
%This ensures Condition \eqref{eq:Cmu} for ${\cal M}_{\alpha,\beta}$ with $p=1<\beta<\alpha=q$.
Moreover for $\alpha=2$, $Z_2(\mu_t)$ is a fractional Brownian motion of Hurst index $H=(3-\beta)/2\in(1/2,1)$ (note that, for $a>0$, $\mu_{at}(A)=a\mu_t(a^{-1}A)$). 
}
\end{Rem}

\begin{Rem}
{\rm
When $\alpha=2$, Theorem \ref{theo:stable} exhibits a Gaussian limit field and generalizes Theorem 2.1 in \cite{BEK} with random weights. 
Indeed, in this case, we have (up to some multiplicative constant) $Z_2=W_\beta$.  
}
\end{Rem}

\begin{Rem}
{\rm A natural complementary result to be investigated is the tightness of $M_\rho$ after normalization and centering
which would allow to turn f.d.d. convergences into weak functional convergences. 
In dimension $d=1$, only partial tightness results are available for the processes studied in \cite{KT}, \cite{MRRS} 
(see Section 4 on "continuous flow reward model" in \cite{KT} and the remarks of Th. 1, Th. 2 and Th. 3 in \cite{MRRS}).
In the case of generalized random fields, tightness issue is more difficult to tackle due to the lack of tractable tightness criterion. 
}
\end{Rem}
%%%%%%%%%%%%%%%%%%%%%%%%%%%%%%%%%%%%%%%%%%%%%%%%%%%%%%%%%%%%%%%%%%%%%%%%%%%%%%%%%%%%%%%%%%%%%%%%%%%%%%%%%%%%%%%%%%

\subsubsection{Poissonian regime}
\label{sec:poisson}
In this section, we investigate the behavior of $M$ under the scaling $\rho^\beta\lambda(\rho)\to a \in (0,\infty)$. 
In this case, the limiting field is given by a compensated Poisson integral and we refer to \cite{Kallenberg} for a general description of Poisson integral. 
We have:
\begin{theo}
\label{theo:poisson}
Suppose $\lambda(\rho)\rho^\beta \to a^{d-\beta}$ when $\rho\to 0^{-\epsilon}$ for some $a>0$. 
We have 
$$
M_\rho(\mu) -\ee[M_\rho(\mu)] \stackrel{{\cal M}_{\alpha,\beta}}{\lto} J(\mu_a),\quad \rho\to0^{-\epsilon}
$$
where $\mu_a$ is the dilatation of $\mu$ and $J$ is the compensated Poisson integral
\begin{equation}
\label{eq:J}
J(\mu)=\int_{\real^d\times\real\times\real^+} m\mu(B(x,r)) \widetilde N_\beta(dx,dr,dm)
\end{equation}
with respect to the compensated Poisson random measure $\widetilde N_\beta $ with intensity given by $C_\beta r^{-\beta-1}dxdrG(dm)$.
\end{theo}

\noindent
Note that the Poisson integral in \eqref{eq:J} above is well defined since  
\begin{equation}
\label{eq:PoissonJ}
\int_{\real^d\times\real\times\real^+}\left(|m\mu(B(x,r))| \wedge (m\mu(B(x,r)))^2\right)%C_\beta a 
r^{-\beta-1}dxdrG(dm)<+\infty
\end{equation}
see Section \ref{sec:Poisson-existence}. 
As  the stable field $Z_\alpha$, the Poisson field $J$ enjoys similar properties.
However, note that in contrast to $Z_\alpha$, $J$ is not self-similar but (and similarly to \cite{BEK}, see also \cite{Kaj05})
$J$ satisfies an aggregate similarity property. 
\begin{prop}
\label{prop:properties_poisson}
The field $J$ is stationary and isotropic on ${\cal M}_{\alpha,\beta}$. 
Moreover, $J$ is aggregate similar, {\it viz.}
$\forall \mu\in{\cal M}_{\alpha,\beta}, \forall m\geq 1,$
\begin{equation}
\label{eq:aggregate}
J(\mu_{a_m})\stackrel{fdd}{=}\sum_{i=1}^m J^i(\mu)
\end{equation}
where $J^i$, $1\leq i\leq m$, are independent copies of $J$ and $a_m=m^{1/(d-\beta)}$.
\end{prop}
The proof of this proposition follows from straightforward computation and will be omitted. 
A comparison of the limiting procedures in Theorem \ref{theo:stable} where $\lambda(\rho)\rho^\beta\to+\infty$ 
and in Theorem \ref{theo:poisson} where $\lambda(\rho)\rho^\beta\to a^{d-\beta}$ suggests that when $a^{d-\beta}\to+\infty$, 
we can recover $Z_\alpha$ from $J$. 
This is true and precisely stated in the following proposition: 
\begin{prop}
\label{prop:ZJ}
When $a^{d-\beta}\to +\infty$, we have 
$\frac 1{a^{(d-\beta)/\alpha}}J(\mu_a)\stackrel{{\cal M}_{\alpha,\beta}}{\lto} Z_\alpha(\mu)$.
\end{prop}

\begin{Rem}
\label{rem:KT_Poisson}
{\rm As in Remark \ref{rem:KT_stable_dep}, when $d=1$ and $\mu_t=|\cdot\cap(0,t)|$, 
the field $J(\mu_t)$ coincides with the {\it intermediate Telecom process} obtained in the intermediate connection rate for the "continuous flow reward model"
in \cite[Th. 1]{KT}, see also Remark \ref{rem:KT1} above.
}
\end{Rem}

\begin{Rem}
{\rm
When $\alpha=2$, Theorem \ref{theo:poisson} generalizes Theorem 2.5 in \cite{BEK} with random weights. 
The field $J$ recovers $J_\beta$ in \cite{BEK} when the random weights in our model are constant. 
Otherwise the law of $J$ depends on the law $G$ of the weight. 
}
\end{Rem}

%%%%%%%%%%%%%%%%%%%%%%%%%%%%%%%%%%%%%%%%%%%%%%%%%%%%%%%%%%%%%%%%%%%%%%%%%%%%%%%%%%%%%%%%%%%%%%%%%%%%%%%%%%%%%%%%%%

\subsubsection{Stable regime with independence for small radius}
\label{sec:stable2}
In this section, we investigate the behavior of $M$ under the scaling $\rho^\beta\lambda(\rho)\to 0$, 
but we restrict to the case $d<\beta<\alpha d$, i.e. $\epsilon=-1$ and $\rho\to 0$, see Section \ref{sec:extra} for $\beta>\alpha d$.  
The case $m\equiv 1$ is considered in Theorem 2 $iii)$ of \cite{KLNS} and we extend here the results and proofs to the case when the weights are random and belong to the normal domain of attraction of a stable distribution.
In comparison to the case $\beta<d$, the tails of the law of the radius are lighter and thus the radius considered are small.  
We show that the asymptotic behavior is given again by a stable field but with index $\gamma=\beta/d$ and defined on $\real^d$.
Moreover in contrast to the stable field $Z_\alpha$ of Section \ref{sec:stable1}, this new field exhibits independence. 
\begin{theo}
\label{theo:indep-stable}
Let $d<\beta<\alpha d$ and suppose that $\lambda(\rho)\rightarrow +\infty$ and $\lambda(\rho)\rho^\beta \to 0$ as $\rho\to 0$. 
Then with $n_2(\rho):=\lambda(\rho)^{d/\beta}\rho^d$ and $\gamma=\beta/d\in(1,\alpha)$, we have 
$$
\frac{M_\rho(\cdot) -\ee[M_\rho(\cdot)]}{n_2(\rho)} \stackrel{L^1(\real^d)\cap L^\alpha(\real^d)}{-\!\!\!-\!\!\!-\!\!\!-\!\!\!-\!\!\!-\!\!\!-\!\!\!-\!\!\!-\!\!\!\lto} \widetilde Z_\gamma(\cdot)
$$
where, for $\mu(dx)=\phi(x)dx$,
$\widetilde Z_\gamma(\mu)=\int_{\real^d} \phi(x) \widetilde M_\gamma(dx)$ 
is a stable integral with respect to the $\gamma$-stable measure $M_\gamma$ with control measure $\sigma_\gamma^\gamma dx$ for 
$$
\sigma_\gamma^\gamma=\frac{c_d^\gamma C_\beta}d \int_{\real_+} \frac{1-\cos(r)}{r^{1+\gamma}}dr \int_\real |m|^\gamma G(dm)
$$
and with constant skewness function equals to
\begin{equation}\label{eq:bgamma}
b_\gamma=-\frac{\int_{\real} \varepsilon(m)|m|^\gamma G(dm)}{\int_{\real} |m|^\gamma G(dm)}.
\end{equation}
\end{theo}
Note that the integrals above are well defined when $d<\beta<\alpha d$ (see Lemma \ref{lemme:fcstable} below). 
The limiting field $\widetilde Z_\gamma$ enjoys similar properties as $Z_\alpha$ and~$J$:
\begin{prop}
The field $\widetilde Z_\gamma$ is stationary, isotropic and self-similar with index $(d-\beta)/\gamma$. 
\end{prop}

\begin{Rem}
\label{rem:KT_stable_indep}
{\rm
As in Remarks \ref{rem:KT_stable_dep} and \ref{rem:KT_Poisson}, 
when $d=1$ and $\phi_t=\ind_{(0,t)}$, the field $\widetilde Z_\gamma(\phi_t)$ coincides with the process obtained in the slow connection rate for the "continuous flow reward model" 
in \cite[Th. 3]{KT}, see also Remark \ref{rem:KT1} above. In this particular case, $\widetilde Z_\gamma(\phi_t)$ is a $\gamma$-stable Lévy process. 
}
\end{Rem}

%%%%%%%%%%%%%%%%%%%%%%%%%%%%%%%%%%%%%%%%%%%%%%%%%%%%%%%%%%%%%%%%%%%%%%%%%%%%%%%%%%%%%%%%%%%%%%%%%%%%%%%%%%%%%%%%%%

\subsubsection{Stable regime with independence for very small radius}
\label{sec:extra}
When the tails of the radii are lighter than that in Section \ref{sec:stable2}, i.e. $\beta>\alpha d$, 
the same stable regime with independence as in Section \ref{sec:stable2} appears but under a different normalization $n_3(\rho):=\lambda(\rho)^{1/\alpha}\rho^d$ and a different stability index $\alpha$.
As previously, since $\beta>\alpha d$, we have $\epsilon=-1$ and the limits are taken when $\rho\to 0$, i.e. the limiting scheme is a zooming-out procedure. 

\begin{theo}
\label{theo:indep-stable2}
Let $\beta>\alpha d$ and suppose that $\lambda(\rho)\to+\infty$ as $\rho\to 0$. 
Let $n_3(\rho):=\lambda(\rho)^{1/\alpha}\rho^d$, then
$$
\frac{M_\rho(\cdot) -\ee[M_\rho(\cdot)]}{n_3(\rho)} \stackrel{L^1(\real^d)\cap L^\alpha(\real^d)}{-\!\!\!-\!\!\!-\!\!\!-\!\!\!-\!\!\!-\!\!\!-\!\!\!-\!\!\!-\!\!\!\lto} \widetilde Z_\alpha(\cdot)
$$
where, for $\mu(dx)=\phi(x)dx$,
$\widetilde Z_\alpha(\mu)=\int_{\real^d} \phi(x) \widetilde M_\alpha(dx)$ 
is a stable integral with respect to the $\alpha$-stable measure $M_\alpha$ with control measure 
$\sigma_\alpha dx$ with $\sigma_\alpha=\sigma c_d \left(\int_{\real^+}r^{\alpha d}F(dr)\right)^{1/\alpha}$
and constant skewness equal to $b$.
\end{theo}

\begin{Rem}
{\rm 
It is worth noting that in both Theorem \ref{theo:indep-stable} and Theorem~\ref{theo:indep-stable2}, 
the stable regime is driven by the parameter $\gamma=(\beta/d)\wedge \alpha$, 
since the normalization is $\lambda(\rho)^{1/\gamma}\rho^d$ and the stability index is $\gamma$. 

Actually, only the asymptotics of the law with the heavier tails contribute to the limit 
while the law with the lighter tails appears only (but globally) as a mere parameter in the limit.
In particular, observe that Theorem~\ref{theo:indep-stable2} applies for any distribution $F$ such that $\int_{\real^+}r^{\alpha d} F(dr)<+\infty$.
}
\end{Rem}

\begin{Rem}
{\rm 
When $d=1$ and $\mu_t=|\cdot\cap(0,t)|$, we recover $ii)$ in Theorem 4 of \cite{KT}.  
}
\end{Rem}

%%%%%%%%%%%%%%%%%%%%%%%%%%%%%%%%%%%%%%%%%%%%%%%%%%%%%%%%%%%%%%%%%%%%%%%%%%%%%%%%%%%%%%%%%%%%%%%%%%%%%%%%%%%%%%%%%%

\section{Proof of the results}
\label{sec:proof}

In what follows, note that the linearity of the random functionals $M_\rho$ and  of the stochastic integrals in $W_\alpha$, $\widetilde W_\alpha$, $Z_\alpha$, $J$ and $\widetilde Z_\gamma$,
 together with the Cram\'er-Wold device imply that the convergence of the finite-dimensional distributions of the centered and renormalized version of $M_\rho$ is equivalent to the convergence of the one-dimensional distributions.
To do so, we will explicitly compute the limits of the characteristic functions, denoting $\varphi_X$ for the characteristic function of a random variable $X$. 
Observe that the characteristic function of $n(\rho)^{-1}(M_\rho(\mu)-\ee[M_\rho(\mu)])$ rewrites:
\begin{eqnarray*}
&&\varphi_{n(\rho)^{-1}(M_\rho(\mu)-\ee[M_\rho(\mu)])}(\theta)\\
&=&\exp\left(\int_{\real^d\times\real^+\times\real} \Psi\left(n(\rho)^{-1}\theta m\mu(B(x,r))\right)\lambda(\rho)dx F_\rho(dr)G(dm)\right)
\end{eqnarray*}
where $\Psi(u)=e^{iu}-1-iu$, see \cite{Kallenberg}. 
Integrating first with respect to the probability $G(dm)$, we have 
\begin{eqnarray}
\nonumber
&&\varphi_{n(\rho)^{-1}(M_\rho(\mu)-\ee[M_\rho(\mu)])}(\theta)\\
\label{eq:fcM}
&=&\exp\left(
\int_{\real^d\times\real^+} \lambda(\rho)\Psi_G\left(n(\rho)^{-1}\theta \mu(B(x,r))\right)dx F_\rho(dr)\right)
\end{eqnarray}
where $\Psi_G(u)=\int_\real \Psi(mu)G(dm)$. 
We also recall that the characteristic function of the stable distribution $S_\alpha(\sigma, b,\tau)$ is given by 
$\exp(-\sigma^\alpha |x|^\alpha(1-ib\varepsilon(\theta) \tan(\pi\alpha/2))+i\tau \theta)$. 

%%%%%%%%%%%%%%%%%%%%%%%%%%%%%%%%%%%%%%%%%%%%%%%%%%%%%%%%%%%%%%%%%%%%%%%%%%%%%%%%%%%%%%%%%%%%%%%%%%%%%%%%%%%%%%%%%%

\subsection{Preliminary lemmas}

In this section, we collect some useful lemmas that will be needed in the proof of our limit theorems \ref{theo:stable}, \ref{theo:poisson} and  \ref{theo:indep-stable}.
We recall the following estimate for the characteristic function of distribution in the domain of attraction of a stable law:

\begin{lemme}
\label{lemme:fcstable}
Suppose $X$ is in the domain of attraction of an $\alpha$-stable law $S_\alpha(\sigma,b,0)$ for some $\alpha>1$.
Then 
$$
\varphi_X(\theta)-1-i\theta \ee[X]\sim_0-\sigma^\alpha|\theta|^\alpha (1-i\varepsilon(\theta) \tan(\pi\alpha/2)b). 
$$
Furthermore, there is some $C>0$ such that for any $\theta\in\real$,
$$
\big|\varphi_X(\theta)-1-i\theta \ee[X]\big|\leq C|\theta|^\alpha.
$$
\end{lemme}

The following lemma is a reformulation from lemma 2.4 in \cite{BEK}. 
It shows that in the scaling limit $\rho\to 0^{-\epsilon}$, 
the behavior of $F_\rho$ is given by the power tail of $F$. 
This is crucial in several estimates.
\begin{lemme}
\label{lemme:BEK}
Let $F$ be as in \eqref{eq:F1} and $\epsilon=\pm 1$.
Assume that $g$ is a continuous function on $\real^+$ such that for some $0<p<\beta<q$, 
there exists some $C>0$ such that 
\begin{equation}
\label{eq:lemmeBEK}
|g(r)|\leq C(r^p\wedge r^q).
\end{equation}
Assume furthermore that $(g_\rho)_{\rho>0}$ is a family of continuous functions such that 
\begin{equation}
\label{eq:lemmeBEK-2}
\lim_{\rho \to 0^{-\epsilon}} |g(r)-g_\rho(r)|=0 \quad {\rm and}\quad  |g(r)-g_\rho(r)|\leq C(r^p\wedge r^q).
\end{equation}
Then 
$$
\int_{\real^+} g_\rho(r)F_\rho(dr)\sim C_\beta \rho^\beta\int_{\real^+} g(r)r^{-1-\beta}dr 
\quad \mbox{ when } \rho \to 0^{-\epsilon}.
$$
\end{lemme}
In the proof of Theorem \ref{theo:stable} and of Theorem \ref{theo:poisson} below, 
this lemma will be used in the particular case where $g_\rho=g$ and $g$ satisfies condition \eqref{eq:lemmeBEK}. 
Roughly speaking, the proof of Lemma \ref{lemme:BEK} consists in taking the limit in the integral. 
This is authorized by the dominated convergence theorem under \eqref{eq:lemmeBEK} and \eqref{eq:lemmeBEK-2}.
We refer to \cite{BEK} for more details.

\subsection{Proofs of Propositions \ref{prop:subspace1} and \ref{prop:subspace2}}

%Note that ${\cal M}_\beta^1$ is indeed a subset of ${\cal M}_\beta$ since for $\mu\in{\cal M}_\beta^1$:
%\begin{eqnarray*}
%\int_{\real^d\times\real^+} |\mu(B(x,r))|^\alpha r^{-1-\beta} drdx
%&=&\int_{\real^+} \left(\int_{\real^d} |\mu(B(x,r))|^\alpha dx \right) r^{-1-\beta} dr\\
%&\leq&\int_{\real^+} C(r^{p-\beta-1} \wedge r^{q-\beta-1})dr
%\end{eqnarray*}
%and $r^{q-\beta-1}$, resp. $r^{p-\beta-1}$, is integrable in $0$, resp. in $+\infty$. 
\medskip\noindent
{\bf Proof of Proposition \ref{prop:subspace1}.}
\medskip\noindent
{\bf Proof of i).}
If \eqref{eq:Cmu} holds true for $\mu_1$ with $p_1<\beta<q_1$ and for $\mu_2$ with $p_2<\beta<q_2$, 
then \eqref{eq:Cmu} holds true for $\mu_1$ and $\mu_2$ with $p=p_1\vee p_2<\beta$ and $q=q_1\wedge q_2>\beta$ (possibly with a different constant $C$). For all $a_1, a_2\in \real$: 
\begin{eqnarray*}
\int_{\real^d}|(a_1\mu_1+a_2\mu_2)(B(x,r))|^\alpha dx
&=&\|(a_1\mu_1+a_2\mu_2)(B(x,r))\|_\alpha^\alpha\\
&\leq& \big(|a_1|\|\mu_1(B(x,r))\|_\alpha+|a_2|\|\mu_2(B(x,r))\|_\alpha\big)^\alpha\\
&\leq& \big((|a_1|^\alpha C(r^p\wedge r^q))^{1/\alpha}+(|a_2|^\alpha C(r^p\wedge r^q))^{1/\alpha}\big)^{\alpha}\\
&=& C(|a_1|+|a_2|)^\alpha(r^p\wedge r^q).
\end{eqnarray*}
This is \eqref{eq:Cmu} for $a_1\mu_1+a_2\mu_2$. 

\medskip\noindent
{\bf Proof of ii).} 
Since $(\tau_s\mu)(B(x,r))=\mu(B(x-s,r))$, $(\theta\mu)(B(x,r))=\mu(B(\Theta^{-1} x,r))$, $\mu_a(B(x,r))=\mu(B(a^{-1}x,a^{-1}r))$, 
the closeness of ${\cal M}_{\alpha,\beta}$ by translations $\tau_s$, by rotations $\Theta$ and by dilatations $x\mapsto ax$ follow straightforwardly from 
the invariance of the Lebesgue measure by translations, by rotation, and by an immediate change of variable in \eqref{eq:Cmu}.

\medskip\noindent
{\bf Proof of iii).}
Since $|\mu|(\real^d)<+\infty$, for $\mu\in {\cal M}_{\alpha,\beta}$ and $\alpha\leq\alpha'$, we have 
\begin{eqnarray*}
\int_{\real^d} |\mu(B(x,r))|^{\alpha'} dx&=&\int_{\real^d} |\mu(B(x,r))|^{\alpha'-\alpha} |\mu(B(x,r))|^{\alpha}dx\\
&\leq& |\mu|(\real^d)^{\alpha'-\alpha}\int_{\real^d}  |\mu(B(x,r))|^{\alpha}dx\\
&\leq& C (r^p\wedge r^q)
\end{eqnarray*}
which proves $\mu\in {\cal M}_{\alpha',\beta}$.

\medskip\noindent
{\bf Proof of iv).}
We prove that $\mu\in {\cal M}_{\alpha,\beta}$ is diffuse when $\beta>d$.
Indeed, suppose that $\mu$ has an atom $a$, then for small enough $r$, $\gamma(r)\geq |\mu(a)/2|^\alpha c_d r^d$, 
where we recall that $\gamma(r)$ is defined in \eqref{eq:Cmu}.
To see this, let $\varepsilon>0$ be such that $\big||\mu|(B(a,\varepsilon))-|\mu(a)|\big|<|\mu(a)|/2$. Then, for every $r<\varepsilon/2$ and $x\in B(a,r)$, $|\mu(B(x,r))|\geq |\mu(a)|/2$. Integrating on $x\in B(a,r)$, we get $\gamma(r)\geq (|\mu(a)|/2)^{\alpha} c_d r^d$.
%let $\delta>0$ be such that $B(a,\delta)\cap \mbox{Supp }(\mu)=\{a\}$. For $r<\delta/2$,
%\begin{eqnarray}
%\nonumber
%\gamma(r)&=&\int_{\real^d}|\mu(B(x,r))|^\alpha dx\\
%\nonumber&\geq &\int_{B(a,r)}|\mu(B(x,r))|^\alpha dx\\
%\nonumber
%&= &\int_{B(a,r)}|\mu(a)|^\alpha dx\\
%\label{eq:atom}
%&=&c_d |\mu(a)|^\alpha r^d.
%\end{eqnarray}
This is in contradiction with \eqref{eq:Cmu} which rewrites $\gamma(r)\leq Cr^q$ for $q>\beta>d$ when $r$ is small.

\medskip\noindent
{\bf Proof of v).}
We prove that $\mu\in{\cal M}_{\alpha, \beta}$ is centered when $\beta\leq d$. 
We will show that 
\begin{equation}
\label{eq:proofv}
\gamma(r)\geq |\mu(\real^d)/3|^\alpha c_dr^d
\end{equation}
when $r$ is large enough.
This is sufficient to prove $v)$ since \eqref{eq:Cmu} rewrites $\gamma(r)\leq Cr^p$ for $p<\beta<d $ when $r\geq 1$
which is in contradiction with \eqref{eq:proofv} when $\mu(\real^d)\not=0$. 

The bound \eqref{eq:proofv} is obvious if $\mu(\real^d)=0$ and if $\mu(\real^d)\neq 0$, let $M$ be such that $|\mu|(B(0,M)^c)\leq |\mu(\real^d)|/3$. 
Then, for $r\geq M$ and any $x\in B(0,r-M)$, $B(0,M)\subset B(x,r)$ and $|\mu(B(x,r))|\geq |\mu(\real^d)|-|\mu|(B(x,r)^c)\geq 2|\mu(\real^d)|/3$. 
Integrating on $x\in B(0,r-M)$, we obtain $\gamma(r)\geq (2|\mu(\real^d)|/3)^\alpha c_d(r-M)^d$. 
This implies \eqref{eq:proofv}. 

\bigskip\noindent
{\bf Proof of Proposition \ref{prop:subspace2}.}
\medskip\noindent
{\bf Proof of i).}
First, when $d<\beta<\alpha d$ and $\mu(dx)=\phi(x)dx\in L^1(\real^d)\cap L^\alpha(\real^d)$,  we have: 
\begin{eqnarray}
 \nonumber
\int_{\real^d} |\mu(B(x,r))|^\alpha dx
&=&\int_{\real^d} \left|\int_{\real^d}  \ind_{\{\|x-y\|<r\}}\phi(y)dy\right|^\alpha dx\\
 \nonumber
&\leq&\int_{\real^d} \left(\int_{\real^d}  \ind_{\{\|x-y\|<r\}}|\phi(y)|^\alpha dy\right) \left(\int_{\real^d}  \ind_{\{\|x-y\|<r\}}dy\right)^{\alpha-1} dx\\
 \nonumber
&=&(c_dr^d)^{\alpha-1} \int_{\real^d} \int_{\real^d}  \ind_{\{\|x-y\|<r\}}|\phi(y)|^\alpha dy dx\\
 \nonumber
&=&(c_dr^d)^{\alpha-1} \int_{\real^d} \left(\int_{\real^d}  \ind_{\{\|x-y\|<r\}}dx\right)|\phi(y)|^\alpha dy\\
\nonumber
&=&(c_dr^d)^{\alpha} \int_{\real^d} |\phi(y)|^\alpha dy
\label{eq:q}
\end{eqnarray}
where we applied Hölder's inequality with $\alpha>1$. 
%Since 
%$\eta>d^2/\beta>d/\alpha>1/\alpha$ in $iii)$ 
%and $f$ is compactly supported in $ii)$, 
Next, in the same way,
\begin{eqnarray}
 \nonumber
\int_{\real^d} |\mu(B(x,r))|^\alpha dx
&=&\int_{\real^d} \left|\int_{\real^d}  \ind_{\{\|x-y\|<r\}}\phi(y)dy\right|^\alpha dx\\
 \nonumber
&\leq&\int_{\real^d} \left(\int_{\real^d}  \ind_{\{\|x-y\|<r\}}|\phi(y)|dy\right) \left(\int_{\real^d}  |\phi(y)| dy\right)^{\alpha-1} dx\\
\label{eq:q2}
&=& c_dr^d \left(\int_{\real^d} |\phi(y)|dy\right)^\alpha. 
\end{eqnarray}
As a consequence, condition \eqref{eq:Cmu} holds with $p=d<\beta<q=\alpha d$, and $\mu\in{\cal M}_{\alpha,\beta}$.

\medskip\noindent
{\bf Proof of ii).} 
Suppose $d-1<\beta< d$ and $\mu(dx)=\phi(x)dx\in L^1(\real^d)$ is centered.  Using $\mu(\real^d)=0$, we have:
\begin{eqnarray}
& &\int_{\real^d} |\mu(B(x,r))|^\alpha dx\\ 
\nonumber
&=&\int_{\real^d} \left|\int_{\real^d}  (\ind_{\{\|x-y\|<r\}}-\ind_{\{\|x\|<r\}})\phi(y)dy\right|^\alpha dx\\
 \nonumber
&\leq& \int_{\real^d} \left(\int_{\real^d}  |\ind_{\{\|x-y\|<r\}}-\ind_{\{\|x\|<r\}}|^\alpha|\phi(y)| dy\right) \left(\int_{\real^d}  |\phi(y)|dy\right)^{\alpha-1} dx.\\
 \nonumber
\end{eqnarray}
Let $\theta(z)=|B(0,1)\Delta B(z,1)|$ denotes the volume of the symmetric difference of the balls with unit radius centered at $0$ and at $z\in\real^d$.  We have, 
$$
\int_{x\in\real^d} |\ind_{\{\|x-y\|<r\}}-\ind_{\{\|x\|<r\}}| dx =r^d\theta\left(\frac{y}{r} \right).
$$
The function $\theta$ is continuous, upper bounded by $c_d$ and $\theta(z)=O(\|z\|)$ as $z\to 0$. 
As a consequence, the global estimate $|\theta(z)|\leq C\|z\|$ holds true for some $C>0$. This entails
\begin{eqnarray}
\int_{\real^d} |\mu(B(x,r))|^\alpha dx
&=& \left(\int_{\real^d}  |\phi(y)|dy\right)^{\alpha-1}  \int_{\real^d} r^d\theta\left(\frac{y}{r} \right) \phi(y)|dy\\
\nonumber
&\leq&  \left(\int_{\real^d}  |\phi(y)|dy\right)^{\alpha-1} \int_{\real^d} C \|y\| |\phi(y)|dy  \: r^{d-1}\\
\nonumber
&\leq& Cr^{d-1}.
\end{eqnarray}
As a consequence, condition \eqref{eq:Cmu} holds true with $p=d-1<\beta<q=d$ because \eqref{eq:q2} still holds true, and finally $\mu\in{\cal M}_{\alpha,\beta}$.

Alternatively, if $\mu$ has a finite support $\{a_1, \dots, a_p\}$, let $\delta>0$ such that for $1\leq i\leq p$, $B(a_i,\delta)\cap \mbox{Supp }(\mu)=\{a_i\}$. 
For $r<\delta/2$,
\begin{eqnarray}
\nonumber
\gamma(r)&=&\int_{\real^d}|\mu(B(x,r))|^\alpha dx\\
\nonumber
&=&\sum_{i=1}^p\int_{B(a_i,r)}|\mu(B(x,r))|^\alpha dx\\
\nonumber
&= &\sum_{i=1}^p\int_{B(a_i,r)}|\mu(a_i)|^\alpha dx\\
\label{eq:atom}
&=&c_d \sum_{i=1}^p|\mu(a_i)|^\alpha r^d=O(r^d).
\end{eqnarray}
Next, let $M$ be such that $\mu(B(0,M)^c)=0$
and note that $\mu(B(x,r))=0$ when $B(x,r)\cap B(0,M) =\emptyset$ 
or when $B(0,M)\subset B(x,r)$ since $\mu(\real^d)=0$. 
We derive $\mu(B(x,r))=0$ when $\|x\|\leq r-M$ or when $\|x\|\geq M+r$. 
Since $\mu$ is a finite measure, we have 
\begin{eqnarray*}
\gamma(r)&=&\int_{r-M\leq \|x\|\leq r+M} |\mu(B(x,r))|^\alpha dx\\ 
&\leq&c_d\big((r+M)^d-(r-M)^d\big) (|\mu|(\real^d))^\alpha\\ 
&=&O(r^{d-1}), \quad r\to+\infty.
\end{eqnarray*}
Together with \eqref{eq:atom}, this yields condition \eqref{eq:Cmu} with $p=d-1<\beta$ and $q=d>\beta$. 
\CQFD

\medskip\noindent
\begin{Rem}[On the bound for large radii] 
{\rm Note that in order to derive the bound $\gamma(r)\leq r^p$ for $p<\beta$ when $r$ is large, the existence of a density for $\mu$ is not required. We can instead suppose that $\mu$ satisfy some tail condition:  for some $\tilde\eta>d/\alpha$ 
$$ 
|\mu|(B(0,R)^c)=O(R^{-\tilde\eta}) \ {\rm as}\ R\to+\infty.
$$

%when $d<\beta<\alpha d$, it is enough to assume the following bound for the tail of $\mu$: 
%\begin{equation}
%\label{eq:tail}
%|\mu|(B(0,r)^c)\leq \tilde Mr^{-\tilde \eta}
%\end{equation}
%for $r$ large enough and $\tilde\eta>d/\alpha$. 
%The condition on $f$ in iv) entails \eqref{eq:tail} with $\tilde\eta=\eta-d>0$ and $\tilde M=Mc_d(\eta-d)^{-1}$ and here the condition required on $\eta$ is stronger $\eta>d(1+\alpha)/\alpha$; 
%however in condition \eqref{eq:tail}, the existence of a density for $\mu$ is not required. 

%Indeed, if \eqref{eq:tail} holds true for, say, $r>A$, 
%for any $x\in B(0,A+r)^c$, since $B(x,r)\subset B(0,\|x\|-r)^c$, we have 
%$$
%|\mu(B(x,r))| \leq |\mu|(B(x,r))\leq |\mu|(B(0, \|x\|-r)^c)\leq \tilde M(\|x\|-r)^{-\tilde \eta}
%$$  
%Let $\theta=(\beta+d)/(2d)>1$. 
%For $r$ large enough, we have $r^\theta\geq A$ and 
%\begin{eqnarray*}
%& &\int_{\|x\|\geq r^\theta+r} |\mu(B(x,r))|^\alpha dx \nonumber\\
%&\leq& 
%\tilde M^\alpha \int_{\|x\|\geq \tilde A+r}(\|x\|-r)^{-\alpha \tilde\eta} dx \nonumber\\
%&=& 
%\tilde M^\alpha c_d\int_{r^\theta+r}^{+\infty}(u-r)^{-\alpha \tilde\eta} u^{d-1}du \nonumber\\
%&=& \tilde M^\alpha c_d\int_{r^\theta}^{+\infty} u^{-\alpha \tilde\eta} (u+r)^{d-1} du \nonumber\\
%\label{maj1} 
%&=& \tilde M^\alpha c_d \sum_{k=0}^{d-1} {d-1 \choose k}(\alpha \eta-k-1)^{-1}r^{d-1-k+\theta(k+1-\alpha \tilde\eta)}\\ 
%\label{maj11} 
%&=&O\big(r^{\theta(d-\alpha \tilde\eta)}\big)
%\end{eqnarray*}
%where the convergence of the integral above is ensured by the condition $\tilde\eta>d/\alpha$. 
%We conclude as previously in the proof of iv).
}
\end{Rem}

%%%%%%%%%%%%%%%%%%%%%%%%%%%%%%%%%%%%%%%%%%%%%%%%%%%%%%%%%%%%%%%%%%%%%%%%%%%%%%%%%%%%%%%%%%%%%%%%%%%%%%%%%%%%%%%%%%

\subsection{Proof of Theorem \ref{theo:stable}}

The characteristic function of the stable integral $Z_\alpha(\mu)$ is given by
\begin{eqnarray}
\label{eq:fcZ}
&&\varphi_{Z_\alpha(\mu)} (\theta)\\
\nonumber
&=&  \exp\left(-C_\beta\sigma^\alpha\!\!\!\int_{\real^d\times\real^+} \!\!\!\!\!\!\!|\theta \mu(B(x,r))|^\alpha (1-i\varepsilon(\theta\mu(B(x,r)))\tan(\pi\alpha/2)b)r^{-1-\beta} drdx\right). 
\end{eqnarray}
Since the characteristic function of the Poisson integral $n_1(\rho)^{-1}(M_\rho(\mu)-\ee[M_\rho(\mu)])$ is given by \eqref{eq:fcM}, 
%\begin{equation}
%\label{eq:fcMM}
%\varphi_{\frac{M_\rho(\mu)-\ee[M_\rho(\mu)]}{n_1(\rho)}}(\theta)=\exp\left(\int_{\real^d\times\real^+\times\real} \Psi\left(n_1(\rho)^{-1}\theta m\mu(B(x,r))\right)\lambda(\rho)dx F_\rho(dr)G(dm)\right),
%\end{equation}
comparing \eqref{eq:fcZ} and \eqref{eq:fcM}, it is sufficient to show that
\begin{eqnarray}
 \label{eq:Proof1}
  &&\lim_{\rho\to0^{-\epsilon}}\int_{\real^d\times\real^+} \lambda(\rho)\Psi_G\left(n_1(\rho)^{-1}\theta \mu(B(x,r))\right)dx F_\rho(dr)\\ 
 \nonumber
&=&  -C_\beta\sigma^\alpha\int_{\real^d\times\real^+} |\theta \mu(B(x,r))|^\alpha (1-i\varepsilon(\theta\mu(B(x,r)))\tan(\pi\alpha/2)b)r^{-1-\beta} drdx.
\end{eqnarray}
Since $n_1(\rho)=(\lambda(\rho)\rho^\beta)^{1/\alpha}\to+\infty$, Lemma \ref{lemme:fcstable} applies and yields
\begin{eqnarray*}
&&\lambda(\rho)\Psi_G\left(n_1(\rho)^{-1}\theta \mu(B(x,r))\right)\\ 
&\sim& -\sigma^\alpha \rho^{-\beta}|\theta|^\alpha |\mu(B(x,r))|^\alpha (1-i\varepsilon(\theta\mu(B(x,r)))\tan(\pi\alpha/2)b).
\end{eqnarray*}
Since $|\frac{\theta}{n_1(\rho)}\mu(B(x,r))|\leq\frac{\theta}{n_1(\rho)} |\mu|(\real^d)$, this equivalence relation is uniform both in $x$ and $r$ and can be integrated. 
This yields
\begin{eqnarray}
\label{eq:Proof2}
&& \int_{\real^d\times\real^+}\lambda(\rho)\Psi_G\left(n_1(\rho)^{-1}\theta \mu(B(x,r))\right)dx F_\rho(dr) \\
\nonumber
&\sim & -\sigma^\alpha \rho^{-\beta}|\theta|^\alpha \int_{\real^d\times\real^+}|\mu(B(x,r))|^\alpha  (1-i\varepsilon(\theta\mu(B(x,r)))\tan(\pi\alpha/2)b)dx F_\rho(dr). 
\end{eqnarray}
Finally, Lemma~\ref{lemme:BEK} applies with
$$
g(r)=\int_{\real^d} |\mu(B(x,r))|^\alpha (1-i\varepsilon(\theta\mu(B(x,r)))\tan(\pi\alpha/2)b)dx,
$$ 
note that \eqref{eq:Cmu} implies that $g$ satisfies condition \eqref{eq:lemmeBEK}. 
Consequently,
\begin{eqnarray}
\nonumber
&&\int_{\real^d\times\real^+}|\mu(B(x,r))|^\alpha  (1-i\varepsilon(\theta\mu(B(x,r)))\tan(\pi\alpha/2)b)dx F_\rho(dr)\\ 
 \label{eq:Proof3}
&\sim&  C_\beta\rho^\beta \int_{\real^d\times\real^+}|\mu(B(x,r))|^\alpha (1-i\varepsilon(\theta\mu(B(x,r)))\tan(\pi\alpha/2)b) r^{-\beta-1} dx dr.
\end{eqnarray}
Finally, \eqref{eq:Proof2} and \eqref{eq:Proof3} together imply \eqref{eq:Proof1}, 
and as explained at the beginning of Section \ref{sec:proof}, Theorem \ref{theo:stable}.
\CQFD

%%%%%%%%%%%%%%%%%%%%%%%%%%%%%%%%%%%%%%%%%%%%%%%%%%%%%%%%%%%%%%%%%%%%%%%%%%%%%%%%%%%%%%%%%%%%%%%%%%%%%%%%%%%%%%%%%%

\subsection{Proof of Condition \eqref{eq:PoissonJ}}
\label{sec:Poisson-existence}

We prove that Condition \eqref{eq:PoissonJ} for the existence of $J$ is satisfied. 
Note that this condition splits into: 
\begin{equation}
\label{eq:PoissonJ1}
\int_{|m\mu(B(x,r))|\leq 1}(m\mu(B(x,r)))^2 r^{-\beta-1}dxdrG(dm)<+\infty
\end{equation}
and
\begin{equation}
\label{eq:PoissonJ2}
\int_{|m\mu(B(x,r))|\geq 1}|m\mu(B(x,r))| r^{-\beta-1}dxdrG(dm)<+\infty. 
\end{equation}
We shall use the following Lemma for the truncated moments of a distribution in the normal domain attraction of a stable law: 
\begin{lemme}\label{lem:moment}
Let $G$ be in the normal domain attraction of an $\alpha$-stable law for $\alpha>1$. 
There are $C_1, C_2\in (0,+\infty)$ such that for all $x\geq 0$: 
$$
\int_{|m|\geq x} |m| G(dm) \leq C_1 x^{1-\alpha}
\qquad \mbox{ and } \qquad
\int_{-x}^x m^2 G(dm) \leq C_2 x^{2-\alpha}.
$$ 
\end{lemme}
{\bf Proof of Lemma \ref{lem:moment}.}
From \cite[XVII.5]{Feller}, we have $\int_{-x}^x m^2 G(dm)\sim Cx^{2-\alpha}$ when $x\to+\infty$ 
(note that since $G$ is in the {\it normal} domain of attraction, there is no slowly varying function in this estimate).
But since moreover for $x\in [0,1]$ 
\begin{eqnarray*}
\int_{-x}^x m^2G(dm)
%&=&\int_{\real} \int_0^{m^2} du \ind_{|m|\leq |x|} G(dm)\\
%&=&\int_{\real}\int_{\real_+}\ind_{u\leq m^2\leq x^2} G(dm)dx\\
%&=&\int_{\real_+}G(m: u\leq m^2\leq x^2)dx\\
&=&\int_0^{x^2}G(m: u\leq m^2\leq x^2)du
\leq x^2\leq x^{2-\alpha}
\end{eqnarray*}
and the second part is proved. 

Next, since $\lim_{x\to 0}\int_{|m|>x}|m|G(dm)=\int_{\real}|m|G(dm)<+\infty$ while $x^{1-\alpha}\to+\infty$, $x\to 0$, 
the first part comes from \cite[Eq. (5.21)]{Feller}: 
$$
\int_{|m|>x}|m|G(dm)\sim \frac{2-\alpha}{\alpha-1}\frac{1}x\int_{-x}^x m^2 G(dm)\sim \frac{2-\alpha}{\alpha-1}x^{1-\alpha}, 
\quad x\to+\infty. 
$$
\CQFD

\noindent 
Now, we prove \eqref{eq:PoissonJ1} and \eqref{eq:PoissonJ2}. First for \eqref{eq:PoissonJ1}, we have: 
\begin{eqnarray*}
&&\int_{|m\mu(B(x,r))|\leq 1}(m\mu(B(x,r)))^2 r^{-\beta-1}dxdrG(dm)\\
&\leq& \int_{\real^d\times\real^+} \left(\int_{-1/|\mu(B(x,r))|}^{+1/|\mu(B(x,r))|}m^2G(dm)\right)\mu(B(x,r))^2r^{-\beta-1}dxdr\\
&\leq& C_2\int_{\real^d\times\real^+} |\mu(B(x,r))|^{\alpha-2}\mu(B(x,r))^2r^{-\beta-1}dxdr\\
&\leq& C_2\int_{\real^d\times\real^+} |\mu(B(x,r))|^\alpha r^{-\beta-1}dxdr
\end{eqnarray*}
which is finite when $\mu\in{\cal M}_{\alpha,\beta}$ (see Prop. \ref{prop:subspace1}-$i)$). 
Next for \eqref{eq:PoissonJ2}, we have: 
\begin{eqnarray*}
&&\int_{|m\mu(B(x,r))|\geq 1}|m\mu(B(x,r))| r^{-\beta-1}dxdrG(dm)\\
&\leq& \int_{\real^d\times\real^+} \left(\int_{|m|>1/|\mu(B(x,r))|}|m|G(dm)\right)|\mu(B(x,r))|r^{-\beta-1}dxdr\\
&\leq& C_1\int_{\real^d\times\real^+} |\mu(B(x,r))|^{\alpha-1}|\mu(B(x,r))|r^{-\beta-1}dxdr\\
&\leq& C_1\int_{\real^d\times\real^+} |\mu(B(x,r))|^\alpha r^{-\beta-1}dxdr
\end{eqnarray*}
which, again, is finite when $\mu\in{\cal M}_{\alpha, \beta}$. 
\CQFD

%%%%%%%%%%%%%%%%%%%%%%%%%%%%%%%%%%%%%%%%%%%%%%%%%%%%%%%%%%%%%%%%%%%%%%%%%%%%%%%%%%%%%%%%%%%%%%%%%%%%%%%%%%%%%%%%%%

\subsection{Proof of Theorem \ref{theo:poisson}}

As in the proof of Theorem \ref{theo:stable}, it is enough to consider convergence of one-dimensional marginals.
The characteristic function of the Poisson integral $M_\rho(\mu)-\ee[M_\rho(\mu)]$ is given by \eqref{eq:fcM}  
and that of the generalized random field $J(\mu)$ is given by 
\begin{eqnarray*}
\varphi_{J(\mu_a)} (\theta)
&=& \exp\left(\int_{\real^d\times\real^+\times\real} \Psi(\theta m\mu(B(a^{-1}x,a^{-1}r))) C_\beta r^{-1-\beta}  drdxG(dm)\right)\\
&=& \exp\left(\int_{\real^d\times\real^+} \Psi_G(\theta\mu(B(x,r))) C_\beta a^{d-\beta} r^{-1-\beta}  drdx\right).
\end{eqnarray*}
%We are thus left to prove that under the asymptotic ??,
%\begin{eqnarray}
%&&\lim_{\rho\to0^{-\epsilon}}\int_{\real^d\times\real^+} \lambda(\rho)\Psi_G\left(\theta \mu(B(x,r))\right)dx F_\rho(dr) \nonumber\\
%&=&  \int_{\real^d\times\real^+} \Psi_G(\theta\mu(B(x,r))) C_\beta a r^{-1-\beta}  drdx \label{eq:Proof4}. 
%\end{eqnarray}
From Lemma~\ref{lemme:fcstable}, $|\Psi_G(\theta\mu(B(x,r)))|\leq C|\theta|^\alpha |\mu(B(x,r))|^\alpha$ for some $C>0$, 
so that condition \eqref{eq:lemmeBEK} for $g(r)=\int_{\real^d} \Psi_G(\mu(B(x,r))) dx$ is given again by \eqref{eq:Cmu} when $\mu\in{\cal M}_{\alpha,\beta}$.
Thus, Lemma~\ref{lemme:BEK} applies and together with $\lim_{\rho\to 0^{-\epsilon}}\lambda(\rho)\rho^\beta=a^{d-\beta}$
entail
\begin{eqnarray*}
&&\lim_{\rho\to 0^{-\epsilon}}\int_{\real^d\times\real^+} \Psi_G\left(\theta \mu(B(x,r))\right)dx \lambda(\rho)F_\rho(dr)\\
&=& C_\beta a^{d-\beta} \int_{\real^d\times\real^+} \Psi_G\left(\theta \mu(B(x,r))\right)r^{-\beta-1}dr dx.
\end{eqnarray*}

Since one-dimensional convergence is enough, this achieves the proof of Theorem \ref{theo:poisson}.
\CQFD

%%%%%%%%%%%%%%%%%%%%%%%%%%%%%%%%%%%%%%%%%%%%%%%%%%%%%%%%%%%%%%%%%%%%%%%%%%%%%%%%%%%%%%%%%%%%%%%%%%%%%%%%%%%%%%%%%%

\subsection{Proof of Proposition \ref{prop:ZJ}}
We consider the subsequence $a_m=m^{1/(d-\beta)}$. 
From the aggregate-similarity of the field $J$ (see \eqref{eq:aggregate} in Proposition \ref{prop:properties_poisson}), we have: 
$$
\frac 1{a_m^{(d-\beta)/\alpha}} J(\mu_{a_m})\stackrel{fdd}{=}
\frac 1{m^{1/\alpha}}\sum_{i=1}^{m}J^i(\mu)
$$
for independent copies $J^i$, $1\leq i\leq m$, of $J$. 
But 
\begin{eqnarray*}
\varphi_{m^{-1/\alpha}\sum_{i=1}^{m}J^i(\mu)}(\theta)
&=&\big(\varphi_{J(\mu)}(m^{-1/\alpha}\theta)\big)^m\\
&=&\exp\left(m\int_{\real^d\times\real^+}\Psi_G(m^{-1/\alpha}\theta\mu(B(x,r)))C_\beta r^{-1-\beta}drdx\right),
\end{eqnarray*}
and from Lemma \ref{lemme:fcstable}, 
$$
\Psi_G(m^{-1/\alpha}\theta\mu(B(x,r)))\sim \sigma^\alpha|\theta|^\alpha |\mu(B(x,r))|^\alpha (1-i\epsilon(\theta\mu(B(x,r)) \tan(\pi\alpha/2)b).
$$
The relation above is uniform both in $x$ and $r$ and it is thus integrable with respect to $drdx$. This  yields 
\begin{eqnarray*}
&&\lim_{m\to+\infty}\varphi_{m^{-1/\alpha}\sum_{i=1}^{m}J^i(\mu)}(\theta)\\
&=&\exp\left(C_\beta \sigma^\alpha|\theta|^\alpha\int_{\real^d\times\real^+}\!\!\!\!\!\!
|\mu(B(x,r))|^\alpha \left(1-i\epsilon(\theta\mu(B(x,r)) \tan\left(\frac{\pi\alpha}{2}\right)b\right)r^{-1-\beta}drdx\right).
\end{eqnarray*}
A standard argument completes the proof of convergence in distribution along an arbitrary sequences. 
\CQFD

%%%%%%%%%%%%%%%%%%%%%%%%%%%%%%%%%%%%%%%%%%%%%%%%%%%%%%%%%%%%%%%%%%%%%%%%%%%%%%%%%%%%%%%%%%%%%%%%%%%%%%%%%%%%%%%%%%

\subsection{Proof of Theorem \ref{theo:indep-stable}}

We follow the argument in the proof of Theorem 2 in \cite{KLNS}.
Recall that here $d<\beta<\alpha d$ so that $\epsilon=-1$ and the limits are taken when $\rho\to 0$. 
Again, by linearity, using the Cram\'er-Wold device, it is enough to deal with one-dimensional marginals. 
From \eqref{eq:fcM} with a change of variable, the characteristic function rewrites
\begin{eqnarray*}
&&\varphi_{n_2(\rho)^{-1}(M_\rho(\mu)-\ee[M_\rho(\mu)])}(\theta)\\
&=&\exp\left(\int_{\real^d\times\real^+}\Psi_G\left(\theta n_2(\rho)^{-1}\mu(B(x,n_2(\rho)^{1/d}r))\right)\lambda(\rho)dx F_{\rho n_2(\rho)^{-1/d}}(dr)\right).
\end{eqnarray*}
Let $\mu(dz)=\phi(z)dz$ with $\phi\in L^1(\real^d)\cap L^\alpha(\real^d)$, 
then, from Lemma 4 in \cite{KLNS}, as $n_2(\rho)\to 0$, 
$$
n_2(\rho)^{-1}\mu(B(x,n_2(\rho)^{1/d}r)) \to \phi(x)c_d r^d
$$
$dx$ almost everywhere and 
\begin{equation}\label{eq:phiast}
x\mapsto \phi^\ast(x)=\sup_{v>0}\left(c_d^{-1}v^{-d}|\mu|(B(x,v))\right) \in L^\alpha(\real^d).
\end{equation}
As a consequence, 
\begin{eqnarray}
&&\int_{\real^d\times\real^+}\Psi_G\left(\theta n_2(\rho)^{-1}\mu(B(x,n_2(\rho)^{1/d}r))\right)\lambda(\rho)dx F_{\rho n_2(\rho)^{-1/d}}(dr) \nonumber\\
 \label{eq:Proof6}
&\sim & C_\beta \lambda(\rho)\rho^\beta n_2(\rho)^{-\beta/d} \int_{\real^d\times\real^+}\Psi_G\left(\theta \phi(x)c_d r^d\right)r^{-\beta-1}dr dx.
\end{eqnarray}
To see this, apply Lemma \ref{lemme:BEK} to
$$
g(r)=\int_{\real^d}\Psi_G\left(\theta \phi(x)c_d r^d\right)dx
$$
and to
$$
g_\rho(r)= \int_{\real^d}\Psi_G\left(\theta n_2(\rho)^{-1}\mu(B(x,n_2(\rho)^{1/d}r))\right)dx.
$$
Since $|\Psi_G(u)|\leq C(|u|\wedge |u|^\alpha)$, we have
$$
|g(r)|\leq C\min\left(c_d |\theta|\|\phi\|_{L^1}r^d,c_d^\alpha |\theta|^\alpha\|\phi\|_{L^\alpha}^\alpha r^{\alpha d}\right)
$$
so that condition \eqref{eq:lemmeBEK} is satisfied with $p=d$ and $q=\alpha d$.
Furthermore, since $\Psi_G$ is a $K$-Lipschitzian function for some finite $K$, we get
$$
|g(r)-g_\rho(r)|\leq K c_d r^d|\theta|\int_{\real^d}\left|c_d^{-1}r^{-d}n_2(\rho)^{-1}\mu(B(x,n_2(\rho)^{1/d}r))-\phi(x)\right|dx. 
$$
The integrand $\left|c_d^{-1}r^{-d}n_2(\rho)^{-1}\mu(B(x,n_2(\rho)^{1/d}r))-\phi(x)\right|$ converges to zero $dx$ almost everywhere. 
Since its $L^\alpha$-norm is bounded by $\|\phi^\ast\|_{L^\alpha}+\|\phi\|_{L^\alpha}$, it is uniformly integrable 
and as a consequence, 
$$
\lim_{\rho\to 0}\int_{\real^d}\left|c_d^{-1}r^{-d}n_2(\rho)^{-1}\mu(B(x,n_2(\rho)^{1/d}r))-\phi(x)\right|dx = 0. 
$$
On the other hand, since for some $C>0$, $|\Psi_G(v)|\leq C|v|^\alpha$, we obtain
$$
|g(r)-g_\rho(r)|\leq C(\|\phi^\ast\|_{L^\alpha}+\|\phi\|_{L^\alpha})r^{\alpha d}.
$$
Hence, $g_\rho$ satisfy condition \eqref{eq:lemmeBEK-2} with $p=d$ and $q=\alpha d$. 
This proves \eqref{eq:Proof6}.\\

From the definition of $n_2(\rho)$, $\lambda(\rho)\rho^\beta n_2(\rho)^{-\beta/d}=1$.
Furthermore, by splitting the integration over $\real^d$ into $\{x\in\real^d:\theta\phi(x)\geq 0\}$ and $\{x\in\real^d:\theta\phi(x)<0\}$ and performing a change of variable, we have 
\begin{eqnarray*}
 \int_{\real^d\times\real^+}\Psi_G\left(\theta \phi(x)c_d r^d\right)r^{-\beta-1}dr dx
 =D \int_{\real^d}(\theta \phi(x))_+^\gamma dx +\bar D\int_{\real^d}(\theta \phi(x))_-^\gamma dx,
\end{eqnarray*}
where $\bar D$ is the complex conjugate of $D=d^{-1}c_d^\gamma\int_{\real^+} \Psi_G(r)r^{-\gamma-1}dr$. 
We deduce 
$$
\varphi_{n_2(\rho)^{-1}(M_\rho(\mu)-\ee[M_\rho(\mu)])}(\theta)\\
=\exp\left(-\sigma_\phi^\gamma|\theta|^\gamma\left(1+ib_\phi \varepsilon(\theta)\tan\left(\frac{\pi\gamma}2\right)\right)\right)
$$
where
$$
\sigma_\phi^\gamma=\sigma_\gamma^\gamma\int_{\real^d} |\phi(x)|^\gamma dx, \quad
$$
and 
\begin{equation}
\label{eq:b_phi}
b_\phi=
\frac{\int_{\real_+} r^{-1-\gamma}(r-\sin(r))dr}{\tan(\pi\gamma/2))\int_{\real_+} r^{-1-\gamma}(1-\cos(r))dr}
\frac{\int_{\real} \varepsilon(m)|m|^\gamma G(dm)}{\int_{\real} |m|^\gamma G(dm)}
\frac{\int_{\real^d}\varepsilon(\phi(x)) |\phi(x)|^\gamma dx}{\int_{\real^d} |\phi(x)|^\gamma dx}.
\end{equation}
But since for $\gamma\in (1,2)$, 
$$
\int_0^{+\infty} \frac{e^{ixu}-1-ixu}{x^{1+\gamma}} dx
=|u|^\gamma \frac{\Gamma(2-\gamma)}{(1-\gamma)(2-\gamma)}\left(\cos(\pi\gamma/2)-i\varepsilon(u)\sin(\pi\gamma/2)\right)
$$
see Lemma 2 in \cite[XVII.4]{Feller} (with $p=1, q=0$ therein), the first ratio on the right-hand side \eqref{eq:b_phi} is $-1$ 
and we have 
$b_\phi=b_\gamma\frac{\int_{\real^d}\varepsilon(\phi(x)) |\phi(x)|^\gamma dx}{\int_{\real^d} |\phi(x)|^\gamma dx}$
where $b_\gamma$ is given in \eqref{eq:bgamma}. This achieves the proof of Theorem \ref{theo:indep-stable}.
\CQFD

%%%%%%%%%%%%%%%%%%%%%%%%%%%%%%%%%%%%%%%%%%%%%%%%%%%%%%%%%%%%%%%%%%%%

\subsection{Proof of Theorem \ref{theo:indep-stable2}}
The argument uses the same tools as in the proof of Theorem \ref{theo:indep-stable} and we only give here the main lines.
From \eqref{eq:fcM} and a change of variable, the characteristic function rewrites
\begin{eqnarray*}
&&\varphi_{n_3(\rho)^{-1}(M_\rho(\mu)-\ee[M_\rho(\mu)])}(\theta)\\
&=&\exp\left(\int_{\real^d\times\real^+}\Psi_G\left(\theta n_3(\rho)^{-1}\mu(B(x,\rho r))\right)\lambda(\rho)dx F(dr)\right).
\end{eqnarray*}
Let $\mu(dz)=\phi(z)dz$ with $\phi\in L^1(\real^d)\cap L^\alpha(\real^d)$.
Since as $\rho\to 0$
$$
\theta n_3(\rho)^{-1}\mu(B(x,\rho r))\sim \lambda(\rho)^{-1/\alpha} c_dr^d\phi(x)
$$
and $\lambda(\rho)\to+\infty$ 
\begin{eqnarray*}
&&\lim_{\rho\to 0}\lambda(\rho)\Psi_G\left(\theta n_3(\rho)^{-1}\mu(B(x,\rho r))\right)\\
%&\sim& -\lambda(\rho)\sigma^\alpha |\theta|^\alpha \lambda(\rho)^{-1} c_d^\alpha r^{\alpha d} |\phi(x)|^\alpha \big(1-\varepsilon(\theta\phi(x))\tan(\pi\alpha/2)b\big) \\
&=&-\sigma^\alpha c_d^\alpha |\theta|^\alpha r^{\alpha d} |\phi(x)|^\alpha \left(1-\varepsilon(\theta\phi(x))\tan\left(\frac{\pi\alpha}2\right)b\right)
\end{eqnarray*}
$dx$ almost everywhere, and this latter function is integrable with respect to $dxF(dr)$ since $\phi \in L^\alpha(\real^d)$ and
$\int_{\real^+}r^{\alpha d}F(dr)<+\infty$. 
Furthermore, with $\phi^\ast$ given in \eqref{eq:phiast}, we derive the following bound:
%domination is proved thanks to the maximal function $\phi^\ast$ as in the proof of Theorem \ref{theo:indep-stable}: for some $C>0$,
\begin{eqnarray*}
\left|\lambda(\rho)\Psi_G\left(\theta n_3(\rho)^{-1}\mu(B(x,\rho r))\right)\right|
&\leq & \lambda(\rho) C n_3(\rho)^{-\alpha}|\mu(B(x,\rho r))|^\alpha\\
&\leq& C r^{\alpha d} |\phi^\ast(x)|^\alpha. 
\end{eqnarray*}
This upper bound is independent of $\rho$ and integrable with respect to $dxF(dr)$ since $\phi^\ast\in L^\alpha(\real^d)$. 
The dominated convergence theorem yields:
\begin{eqnarray*}
& &\lim_{\rho\to 0}\int_{\real^d\times\real^+}\Psi_G\left(\theta n_3(\rho)^{-1}\mu(B(x,\rho r))\right)\lambda(\rho)dx F(dr)\\
&=& -\sigma^\alpha c_d^\alpha |\theta|^\alpha \int_{\real^+}r^{\alpha d}F(dr) \int_{\real^d}|\phi(x)|^\alpha (1-\varepsilon(\theta\phi(x))\tan(\pi\alpha/2)b)dx.
\end{eqnarray*}
This proves Theorem \ref{theo:indep-stable2}.
\CQFD

\bigskip
\noindent
{\bf Acknowledgment.} 
The first author thanks the LMA, University of Poitiers, for the kind hospitality during the period of preparation of this work  
and the CNRS for its financial support. Both authors thank an anonymous referee for valuable comments on a preliminary version.

%%%%%%%%%%%%%%%%%%%%%%%%%%%%%%%%%%%%%%%%%%%%%%%%%%%%%%%%%%%%%%%%%%%%%%%%%%%%%%%%%%%%%%%%%%%
%%%%%%%%%%%%%%%%%%%%%%%%%%%%%%%%%%%%%%%%%%%%%%%%%%%%%%%%%%%%%%%%%%%%%%%%%%%%%%%%%%%%%%%%%%%

\end{document}